\documentclass{amsart}
\usepackage{amsfonts}

\newtheorem{theorem}{Theorem}
\theoremstyle{plain}

\newtheorem{proposition}{Proposition}
\newtheorem{remark}{Remark}

\numberwithin{equation}{section}
\input{tcilatex}

\begin{document}
\title[Reverses of the Triangle Inequality]{Reverses of the Triangle
Inequality in Inner Product Spaces}
\author{Sever S. Dragomir}
\address{School of Computer Science and Mathematics\\
Victoria University of Technology\\
PO Box 14428, MCMC 8001\\
Victoria, Australia.}
\email{sever@csm.vu.edu.au}
\urladdr{http://rgmia.vu.edu.au/SSDragomirWeb.html}
\date{April 16, 2004.}
\subjclass[2000]{Primary 46C05; Secondary 26D15.}
\keywords{Triangle inequality, Diaz-Metcalf Inequality, Inner products.}

\begin{abstract}
Some new reverses  for the generalised triangle inequality in inner product
spaces and applications are given. Applications in connection to the Schwarz
inequality are provided as well.
\end{abstract}

\maketitle

\section{Introduction}

In 1966, J.B. Diaz and F.T. Metcalf \cite{DM} proved the following reverse
of the triangle inequality:

\begin{theorem}
\label{ta}Let $a$ be a unit vector in the inner product space $\left(
H;\left\langle \cdot ,\cdot \right\rangle \right) $ over the real or complex
number field $\mathbb{K}$. Suppose that the vectors $x_{i}\in H\backslash
\left\{ 0\right\} ,$ $i\in \left\{ 1,\dots ,n\right\} $ satisfy%
\begin{equation}
0\leq r\leq \frac{\func{Re}\left\langle x_{i},a\right\rangle }{\left\Vert
x_{i}\right\Vert },\ \ \ \ \ i\in \left\{ 1,\dots ,n\right\} .  \label{1.1}
\end{equation}%
Then%
\begin{equation}
r\sum_{i=1}^{n}\left\Vert x_{i}\right\Vert \leq \left\Vert
\sum_{i=1}^{n}x_{i}\right\Vert ,  \label{1.2}
\end{equation}%
where equality holds if and only if%
\begin{equation}
\sum_{i=1}^{n}x_{i}=r\left( \sum_{i=1}^{n}\left\Vert x_{i}\right\Vert
\right) a.  \label{1.3}
\end{equation}
\end{theorem}

A generalisation of this result for orthornormal families is incorporated in
the following result \cite{DM}.

\begin{theorem}
\label{tb}Let $a_{1},\dots ,a_{n}$ be orthornormal vectors in $H.$ Suppose
the vectors $x_{1},\dots ,x_{n}\in H\backslash \left\{ 0\right\} $ satisfy%
\begin{equation}
0\leq r_{k}\leq \frac{\func{Re}\left\langle x_{i},a_{k}\right\rangle }{%
\left\Vert x_{i}\right\Vert },\ \ \ \ \ i\in \left\{ 1,\dots ,n\right\} ,\
k\in \left\{ 1,\dots ,m\right\} .  \label{1.4}
\end{equation}%
Then%
\begin{equation}
\left( \sum_{k=1}^{m}r_{k}^{2}\right) ^{\frac{1}{2}}\sum_{i=1}^{n}\left\Vert
x_{i}\right\Vert \leq \left\Vert \sum_{i=1}^{n}x_{i}\right\Vert ,
\label{1.5}
\end{equation}%
where equality holds if and only if%
\begin{equation}
\sum_{i=1}^{n}x_{i}=\left( \sum_{i=1}^{n}\left\Vert x_{i}\right\Vert \right)
\sum_{k=1}^{m}r_{k}a_{k}.  \label{1.6}
\end{equation}
\end{theorem}

Similar results valid for semi-inner products may be found in \cite{K} and 
\cite{M}.

For other inequalities related to the triangle inequality, see Chapter XVII
of the book \cite{MPF} and the references therein.

The main aim of this paper is to point out new reverses of the generalised
triangle inequality which naturally complement the above results due to Diaz
and Metcalf. New reverses for the Schwarz inequality are provided.
Applications for vector-valued integrals in Hilbert spaces are pointed out
as well.

\section{Some Inequalities of Diaz-Metcalf Type}

The following result with a natural geometrical meaning holds:

\begin{theorem}
\label{t2.1}Let $a$ be a unit vector in the inner product space $\left(
H;\left\langle \cdot ,\cdot \right\rangle \right) $ and $\rho \in \left(
0,1\right) .$ If $x_{i}\in H,$ $i\in \left\{ 1,\dots ,n\right\} $ are such
that%
\begin{equation}
\left\Vert x_{i}-a\right\Vert \leq \rho \text{ \ for each \ }i\in \left\{
1,\dots ,n\right\} ,  \label{2.1}
\end{equation}%
then we have the inequality%
\begin{equation}
\sqrt{1-\rho ^{2}}\sum_{i=1}^{n}\left\Vert x_{i}\right\Vert \leq \left\Vert
\sum_{i=1}^{n}x_{i}\right\Vert ,  \label{2.2}
\end{equation}%
with equality if and only if%
\begin{equation}
\sum_{i=1}^{n}x_{i}=\sqrt{1-\rho ^{2}}\left( \sum_{i=1}^{n}\left\Vert
x_{i}\right\Vert \right) a.  \label{2.3}
\end{equation}
\end{theorem}

\begin{proof}
From (\ref{2.1}) we have%
\begin{equation*}
\left\Vert x_{i}\right\Vert ^{2}-2\func{Re}\left\langle x_{i},a\right\rangle
+1\leq \rho ^{2},
\end{equation*}%
giving%
\begin{equation}
\left\Vert x_{i}\right\Vert ^{2}+1-\rho ^{2}\leq 2\func{Re}\left\langle
x_{i},a\right\rangle ,  \label{2.4}
\end{equation}%
for each $i\in \left\{ 1,\dots ,n\right\} .$

Dividing by $\sqrt{1-\rho ^{2}}>0,$ we deduce%
\begin{equation}
\frac{\left\Vert x_{i}\right\Vert ^{2}}{\sqrt{1-\rho ^{2}}}+\sqrt{1-\rho ^{2}%
}\leq \frac{2\func{Re}\left\langle x_{i},a\right\rangle }{\sqrt{1-\rho ^{2}}}%
,  \label{2.5}
\end{equation}%
for each $i\in \left\{ 1,\dots ,n\right\} .$

On the other hand, by the elementary inequality%
\begin{equation}
\frac{p}{\alpha }+q\alpha \geq 2\sqrt{pq},\ \ \ p,q\geq 0,\ \alpha >0
\label{2.5.1}
\end{equation}%
we have%
\begin{equation}
2\left\Vert x_{i}\right\Vert \leq \frac{\left\Vert x_{i}\right\Vert ^{2}}{%
\sqrt{1-\rho ^{2}}}+\sqrt{1-\rho ^{2}}  \label{2.6}
\end{equation}%
and thus, by (\ref{2.5}) and (\ref{2.6}), we deduce%
\begin{equation*}
\frac{\func{Re}\left\langle x_{i},a\right\rangle }{\left\Vert
x_{i}\right\Vert }\geq \sqrt{1-\rho ^{2}},
\end{equation*}%
for each $i\in \left\{ 1,\dots ,n\right\} .$ Applying Theorem \ref{ta} for $%
r=\sqrt{1-\rho ^{2}},$ we deduce the desired inequality (\ref{2.2}).
\end{proof}

In a similar manner to the one used in the proof of Theorem \ref{t2.1} and
by the use of the Diaz-Metcalf inequality incorporated in Theorem \ref{tb},
we can also prove the following result:

\begin{theorem}
\label{t2.2}Let $a_{1},\dots ,a_{n}$ be orthornormal vectors in $H.$ Suppose
the vectors $x_{1},\dots ,x_{n}\in H\backslash \left\{ 0\right\} $ satisfy%
\begin{equation}
\left\Vert x_{i}-a_{k}\right\Vert \leq \rho _{k}\text{ \ for each}\ \ i\in
\left\{ 1,\dots ,n\right\} ,\ k\in \left\{ 1,\dots ,m\right\} ,  \label{2.11}
\end{equation}%
where $\rho _{k}\in \left( 0,1\right) .$ Then we have the following reverse
of the triangle inequality%
\begin{equation}
\left( m-\sum_{k=1}^{m}\rho _{k}^{2}\right) ^{1/2}\sum_{i=1}^{n}\left\Vert
x_{i}\right\Vert \leq \left\Vert \sum_{i=1}^{n}x_{i}\right\Vert .
\label{2.12}
\end{equation}%
The equality holds in (\ref{2.12}) if and only if%
\begin{equation}
\sum_{i=1}^{n}x_{i}=\left( \sum_{i=1}^{n}\left\Vert x_{i}\right\Vert \right)
\sum_{k=1}^{m}\left( 1-\rho _{k}^{2}\right) ^{1/2}a_{k}.  \label{2.13}
\end{equation}
\end{theorem}

The following result with a different geometrical meaning may be stated as
well:

\begin{theorem}
\label{t2.3}Let $a$ be a unit vector in the inner product space $\left(
H;\left\langle \cdot ,\cdot \right\rangle \right) $ and $M\geq m>0.$ If $%
x_{i}\in H,$ $i\in \left\{ 1,\dots ,n\right\} $ are such that either%
\begin{equation}
\func{Re}\left\langle Ma-x_{i},x_{i}-ma\right\rangle \geq 0  \label{2.14}
\end{equation}%
or, equivalently,%
\begin{equation}
\left\Vert x_{i}-\frac{M+m}{2}\cdot a\right\Vert \leq \frac{1}{2}\left(
M-m\right)  \label{2.15}
\end{equation}%
holds for each $i\in \left\{ 1,\dots ,n\right\} ,$ then we have the
inequality%
\begin{equation}
\frac{2\sqrt{mM}}{m+M}\sum_{i=1}^{n}\left\Vert x_{i}\right\Vert \leq
\left\Vert \sum_{i=1}^{n}x_{i}\right\Vert ,  \label{2.16}
\end{equation}%
or, equivalently, 
\begin{equation}
\left( 0\leq \right) \sum_{i=1}^{n}\left\Vert x_{i}\right\Vert -\left\Vert
\sum_{i=1}^{n}x_{i}\right\Vert \leq \frac{\left( \sqrt{M}-\sqrt{m}\right)
^{2}}{2\sqrt{mM}}\left\Vert \sum_{i=1}^{n}x_{i}\right\Vert .  \label{2.17}
\end{equation}%
The equality holds in (\ref{2.16}) (or in (\ref{2.17})) if and only if%
\begin{equation}
\sum_{i=1}^{n}x_{i}=\frac{2\sqrt{mM}}{m+M}\left( \sum_{i=1}^{n}\left\Vert
x_{i}\right\Vert \right) a.  \label{2.19}
\end{equation}
\end{theorem}

\begin{proof}
Firstly, we remark that if $x,z,Z\in H,$ then the following statements are
equivalent

\begin{enumerate}
\item[(i)] $\func{Re}\left\langle Z-x,x-z\right\rangle \geq 0$

and

\item[(ii)] $\left\Vert x-\frac{Z+z}{2}\right\Vert \leq \frac{1}{2}%
\left\Vert Z-z\right\Vert .$
\end{enumerate}

Using this fact, one may simply realize that (\ref{2.14}) and (\ref{2.15})
are equivalent.

Now, from (\ref{2.14}), we get%
\begin{equation*}
\left\Vert x_{i}\right\Vert ^{2}+mM\leq \left( M+m\right) \func{Re}%
\left\langle x_{i},a\right\rangle ,
\end{equation*}%
for any $i\in \left\{ 1,\dots ,n\right\} .$ Dividing this inequality by $%
\sqrt{mM}>0,$ we deduce the following inequality that will be used in the
sequel%
\begin{equation}
\frac{\left\Vert x_{i}\right\Vert ^{2}}{\sqrt{mM}}+\sqrt{mM}\leq \frac{M+m}{%
\sqrt{mM}}\func{Re}\left\langle x_{i},a\right\rangle ,  \label{2.20}
\end{equation}%
for each $i\in \left\{ 1,\dots ,n\right\} .$

Using the inequality (\ref{2.5.1}) from Theorem \ref{t2.1}, we also have%
\begin{equation}
2\left\Vert x_{i}\right\Vert \leq \frac{\left\Vert x_{i}\right\Vert ^{2}}{%
\sqrt{mM}}+\sqrt{mM},  \label{2.21}
\end{equation}%
for each $i\in \left\{ 1,\dots ,n\right\} .$

Utilizing (\ref{2.20}) and (\ref{2.21}), we may conclude with the following
inequality%
\begin{equation*}
\left\Vert x_{i}\right\Vert \leq \frac{M+m}{\sqrt{mM}}\func{Re}\left\langle
x_{i},a\right\rangle ,
\end{equation*}%
which is equivalent to%
\begin{equation}
\frac{2\sqrt{mM}}{m+M}\leq \frac{\func{Re}\left\langle x_{i},a\right\rangle 
}{\left\Vert x_{i}\right\Vert }  \label{2.22}
\end{equation}%
for any $i\in \left\{ 1,\dots ,n\right\} .$

Finally, on applying the Diaz-Metcalf result in Theorem \ref{ta} for $r=%
\frac{2\sqrt{mM}}{m+M}$, we deduce the desired conclusion.

The equivalence between (\ref{2.16}) and (\ref{2.17}) follows by simple
calculation and we omit the details.
\end{proof}

Finally, by the use of Theorem \ref{tb} and a similar technique to that
employed in the proof of Theorem \ref{t2.3}, we may state the following
result:

\begin{theorem}
\label{t2.4}Let $a_{1},\dots ,a_{n}$ be orthornormal vectors in $H.$ Suppose
the vectors $x_{1},\dots ,x_{n}\in H\backslash \left\{ 0\right\} $ satisfy%
\begin{equation}
\func{Re}\left\langle M_{k}a_{k}-x_{i},x_{i}-\mu _{k}a_{k}\right\rangle \geq
0,  \label{2.23}
\end{equation}%
or, equivalently,%
\begin{equation}
\left\Vert x_{i}-\frac{M_{k}+\mu _{k}}{2}a_{k}\right\Vert \leq \frac{1}{2}%
\left( M_{k}-\mu _{k}\right) ,  \label{2.24}
\end{equation}%
for any $i\in \left\{ 1,\dots ,n\right\} $ and $k\in \left\{ 1,\dots
,m\right\} ,$ where $M_{k}\geq \mu _{k}>0$ for each $k\in \left\{ 1,\dots
,m\right\} .$

Then we have the inequality%
\begin{equation}
2\left( \sum_{k=1}^{m}\frac{\mu _{k}M_{k}}{\left( \mu _{k}+M_{k}\right) ^{2}}%
\right) ^{\frac{1}{2}}\sum_{i=1}^{n}\left\Vert x_{i}\right\Vert \leq
\left\Vert \sum_{i=1}^{n}x_{i}\right\Vert .  \label{2.25}
\end{equation}%
The equality holds in (\ref{2.25}) iff%
\begin{equation}
\sum_{i=1}^{n}x_{i}=2\left( \sum_{i=1}^{n}\left\Vert x_{i}\right\Vert
\right) \sum_{k=1}^{m}\frac{\sqrt{\mu _{k}M_{k}}}{\mu _{k}+M_{k}}a_{k}.
\label{2.26}
\end{equation}
\end{theorem}

\section{Some New Reverses of the Triangle Inequality}

In this section we establish some additive reverses of the generalised
triangle inequality in real or complex inner product spaces.

The following result holds:

\begin{theorem}
\label{t3.1}Let $\left( H;\left\langle \cdot ,\cdot \right\rangle \right) $
be an inner product space over the real or complex number field $\mathbb{K}$
and $e,$ $x_{i}\in H,$ $i\in \left\{ 1,\dots ,n\right\} $ with $\left\Vert
e\right\Vert =1.$ If $k_{i}\geq 0$, $i\in \left\{ 1,\dots ,n\right\} ,$ are
such that%
\begin{equation}
\left\Vert x_{i}\right\Vert -\func{Re}\left\langle e,x_{i}\right\rangle \leq
k_{i}\text{\ \ for each}\ \ i\in \left\{ 1,\dots ,n\right\} ,  \label{3.1}
\end{equation}%
then we have the inequality%
\begin{equation}
\left( 0\leq \right) \sum_{i=1}^{n}\left\Vert x_{i}\right\Vert -\left\Vert
\sum_{i=1}^{n}x_{i}\right\Vert \leq \sum_{i=1}^{n}k_{i}.  \label{3.2}
\end{equation}%
The equality holds in (\ref{3.2}) if and only if%
\begin{equation}
\sum_{i=1}^{n}\left\Vert x_{i}\right\Vert \geq \sum_{i=1}^{n}k_{i}
\label{3.3}
\end{equation}%
and%
\begin{equation}
\sum_{i=1}^{n}x_{i}=\left( \sum_{i=1}^{n}\left\Vert x_{i}\right\Vert
-\sum_{i=1}^{n}k_{i}\right) e.  \label{3.4}
\end{equation}
\end{theorem}

\begin{proof}
If we sum in (\ref{3.1}) over $i$ from 1 to $n,$ then we get%
\begin{equation}
\sum_{i=1}^{n}\left\Vert x_{i}\right\Vert \leq \func{Re}\left\langle
e,\sum_{i=1}^{n}x_{i}\right\rangle +\sum_{i=1}^{n}k_{i}.  \label{3.5}
\end{equation}%
By Schwarz's inequality for $e$ and $\sum_{i=1}^{n}x_{i},$ we have%
\begin{align}
\func{Re}\left\langle e,\sum_{i=1}^{n}x_{i}\right\rangle & \leq \left\vert 
\func{Re}\left\langle e,\sum_{i=1}^{n}x_{i}\right\rangle \right\vert
\label{3.6} \\
& \leq \left\vert \left\langle e,\sum_{i=1}^{n}x_{i}\right\rangle
\right\vert \leq \left\Vert e\right\Vert \left\Vert
\sum_{i=1}^{n}x_{i}\right\Vert =\left\Vert \sum_{i=1}^{n}x_{i}\right\Vert . 
\notag
\end{align}%
Making use of (\ref{3.5}) and (\ref{3.6}), we deduce the desired inequality (%
\ref{3.1}).

If (\ref{3.3}) and (\ref{3.4}) hold, then%
\begin{equation*}
\left\Vert \sum_{i=1}^{n}x_{i}\right\Vert =\left\vert
\sum_{i=1}^{n}\left\Vert x_{i}\right\Vert -\sum_{i=1}^{n}k_{i}\right\vert
\left\Vert e\right\Vert =\sum_{i=1}^{n}\left\Vert x_{i}\right\Vert
-\sum_{i=1}^{n}k_{i},
\end{equation*}%
and the equality in the second part of (\ref{3.2}) holds true.

Conversely, if the equality holds in (\ref{3.2}), then, obviously (\ref{3.3}%
) is valid and we need only to prove (\ref{3.4}).

Now, if the equality holds in (\ref{3.2}) then it must hold in (\ref{3.1})
for each $i\in \left\{ 1,\dots ,n\right\} $ and also must hold in any of the
inequalities in (\ref{3.6}).

It is well known that in Schwarz's inequality $\left\vert \left\langle
u,v\right\rangle \right\vert \leq \left\Vert u\right\Vert \left\Vert
v\right\Vert $ $\left( u,v\in H\right) $ the case of equality holds iff
there exists a $\lambda \in \mathbb{K}$ such that $u=\lambda v.$ We note
that in the weaker inequality $\func{Re}\left\langle u,v\right\rangle \leq
\left\Vert u\right\Vert \left\Vert v\right\Vert $ the case of equality holds
iff $\lambda \geq 0$ and $u=\lambda v.$

Consequently, the equality holds in all inequalities (\ref{3.6})
simultaneously iff there exists a $\mu \geq 0$ with%
\begin{equation}
\mu e=\sum_{i=1}^{n}x_{i}.  \label{3.7}
\end{equation}

If we sum the equalities in (\ref{3.1}) over $i$ from 1 to $n,$ then we
deduce%
\begin{equation}
\sum_{i=1}^{n}\left\Vert x_{i}\right\Vert -\func{Re}\left\langle
e,\sum_{i=1}^{n}x_{i}\right\rangle =\sum_{i=1}^{n}k_{i}.  \label{3.8}
\end{equation}%
Replacing $\sum_{i=1}^{n}\left\Vert x_{i}\right\Vert $ from (\ref{3.7}) into
(\ref{3.8}), we deduce%
\begin{equation*}
\sum_{i=1}^{n}\left\Vert x_{i}\right\Vert -\mu \left\Vert e\right\Vert
^{2}=\sum_{i=1}^{n}k_{i},
\end{equation*}%
from where we get $\mu =\sum_{i=1}^{n}\left\Vert x_{i}\right\Vert
-\sum_{i=1}^{n}k_{i}.$ Using (\ref{3.7}), we deduce (\ref{3.4}) and the
theorem is proved.
\end{proof}

If we turn our attention to the case of orthogonal families, then we may
state the following result as well.

\begin{theorem}
\label{t3.2}Let $\left( H;\left\langle \cdot ,\cdot \right\rangle \right) \ $%
be an inner product space over the real or complex number field $\mathbb{K}$%
, $\left\{ e_{k}\right\} _{k\in \left\{ 1,\dots ,m\right\} }$ a family of
orthonormal vectors in $H,$ $x_{i}\in H,$ $M_{i,k}\geq 0$ for $i\in \left\{
1,\dots ,n\right\} $ and $k\in \left\{ 1,\dots ,m\right\} $ such that%
\begin{equation}
\left\Vert x_{i}\right\Vert -\func{Re}\left\langle e_{k},x_{i}\right\rangle
\leq M_{ik}\text{ \ for each}\ \ i\in \left\{ 1,\dots ,n\right\} ,\ k\in
\left\{ 1,\dots ,m\right\} .  \label{3.9}
\end{equation}%
Then we have the inequality%
\begin{equation}
\sum_{i=1}^{n}\left\Vert x_{i}\right\Vert \leq \frac{1}{\sqrt{m}}\left\Vert
\sum_{i=1}^{n}x_{i}\right\Vert +\frac{1}{m}\sum_{i=1}^{n}%
\sum_{k=1}^{m}M_{ik}.  \label{3.10}
\end{equation}%
The equality holds true in (\ref{3.10}) if and only if 
\begin{equation}
\sum_{i=1}^{n}\left\Vert x_{i}\right\Vert \geq \frac{1}{m}%
\sum_{i=1}^{n}\sum_{k=1}^{m}M_{ik}  \label{3.11}
\end{equation}%
and%
\begin{equation}
\sum_{i=1}^{n}x_{i}=\left( \sum_{i=1}^{n}\left\Vert x_{i}\right\Vert -\frac{1%
}{m}\sum_{i=1}^{n}\sum_{k=1}^{m}M_{ik}\right) \sum_{k=1}^{m}e_{k}.
\label{3.12}
\end{equation}
\end{theorem}

\begin{proof}
If we sum over $i$ from 1 to $n$ in (\ref{3.9}), then we obtain%
\begin{equation*}
\sum_{i=1}^{n}\left\Vert x_{i}\right\Vert \leq \func{Re}\left\langle
e,\sum_{i=1}^{n}x_{i}\right\rangle +\sum_{i=1}^{n}M_{ik},
\end{equation*}%
for each $k\in \left\{ 1,\dots ,m\right\} .$ Summing these inequalities over 
$k$ from 1 to $m,$ we deduce%
\begin{equation}
\sum_{i=1}^{n}\left\Vert x_{i}\right\Vert \leq \frac{1}{m}\func{Re}%
\left\langle \sum_{k=1}^{m}e_{k},\sum_{i=1}^{n}x_{i}\right\rangle +\frac{1}{m%
}\sum_{i=1}^{n}\sum_{k=1}^{m}M_{ik}.  \label{3.13}
\end{equation}%
By Schwarz's inequality for $\sum_{k=1}^{m}e_{k}$ and $\sum_{i=1}^{n}x_{i}$
we have%
\begin{align}
\func{Re}\left\langle \sum_{k=1}^{m}e_{k},\sum_{i=1}^{n}x_{i}\right\rangle &
\leq \left\vert \func{Re}\left\langle
\sum_{k=1}^{m}e_{k},\sum_{i=1}^{n}x_{i}\right\rangle \right\vert 
\label{3.14} \\
& \leq \left\vert \left\langle
\sum_{k=1}^{m}e_{k},\sum_{i=1}^{n}x_{i}\right\rangle \right\vert   \notag \\
& \leq \left\Vert \sum_{k=1}^{m}e_{k}\right\Vert \left\Vert
\sum_{i=1}^{n}x_{i}\right\Vert   \notag \\
& =\sqrt{m}\left\Vert \sum_{i=1}^{n}x_{i}\right\Vert ,  \notag
\end{align}%
since, obviously,%
\begin{equation*}
\left\Vert \sum_{k=1}^{m}e_{k}\right\Vert =\sqrt{\left\Vert
\sum_{k=1}^{m}e_{k}\right\Vert ^{2}}=\sqrt{\sum_{k=1}^{m}\left\Vert
e_{k}\right\Vert ^{2}}=\sqrt{m}.
\end{equation*}%
Making use of (\ref{3.13}) and (\ref{3.14}), we deduce the desired
inequality (\ref{3.10}).

If (\ref{3.11}) and (\ref{3.12}) hold, then%
\begin{align*}
\frac{1}{\sqrt{m}}\left\Vert \sum_{i=1}^{n}x_{i}\right\Vert & =\left\vert
\sum_{i=1}^{n}\left\Vert x_{i}\right\Vert -\frac{1}{m}\sum_{i=1}^{n}%
\sum_{k=1}^{m}M_{ik}\right\vert \left\Vert \sum_{k=1}^{m}e_{k}\right\Vert \\
& =\frac{\sqrt{m}}{\sqrt{m}}\left( \sum_{i=1}^{n}\left\Vert x_{i}\right\Vert
-\frac{1}{m}\sum_{i=1}^{n}\sum_{k=1}^{m}M_{ik}\right) \\
& =\sum_{i=1}^{n}\left\Vert x_{i}\right\Vert -\frac{1}{m}\sum_{i=1}^{n}%
\sum_{k=1}^{m}M_{ik},
\end{align*}%
and the equality in (\ref{3.10}) holds true.

Conversely, if the equality holds in (\ref{3.10}), then, obviously (\ref%
{3.11}) is valid.

Now if the equality holds in (\ref{3.10}), then it must hold in (\ref{3.9})
for each $i\in \left\{ 1,\dots ,n\right\} $ and $k\in \left\{ 1,\dots
,m\right\} $ and also must hold in any of the inequalities in (\ref{3.14}).

It is well known that in Schwarz's inequality $\func{Re}\left\langle
u,v\right\rangle \leq \left\Vert u\right\Vert \left\Vert v\right\Vert ,$ the
equality occurs iff $u=\lambda v$ with $\lambda \geq 0,$ consequently, the
equality holds in all inequalities (\ref{3.14}) simultaneously iff there
exists a $\mu \geq 0$ with%
\begin{equation}
\mu \sum_{k=1}^{m}e_{k}=\sum_{i=1}^{n}x_{i}.  \label{3.15}
\end{equation}%
If we sum the equality in (\ref{3.9}) over $i$ from 1 to $n$ and $k$ from 1
to $m,$ then we deduce%
\begin{equation}
m\sum_{i=1}^{n}\left\Vert x_{i}\right\Vert -\func{Re}\left\langle
\sum_{k=1}^{m}e_{k},\sum_{i=1}^{n}x_{i}\right\rangle
=\sum_{i=1}^{n}\sum_{k=1}^{m}M_{ik}.  \label{3.16}
\end{equation}%
Replacing $\sum_{i=1}^{n}x_{i}$ from (\ref{3.15}) into (\ref{3.16}), we
deduce%
\begin{equation*}
m\sum_{i=1}^{n}\left\Vert x_{i}\right\Vert -\mu \sum_{k=1}^{m}\left\Vert
e_{k}\right\Vert ^{2}=\sum_{i=1}^{n}\sum_{k=1}^{m}M_{ik}
\end{equation*}%
giving%
\begin{equation*}
\mu =\sum_{i=1}^{n}\left\Vert x_{i}\right\Vert -\frac{1}{m}%
\sum_{i=1}^{n}\sum_{k=1}^{m}M_{ik}.
\end{equation*}%
Using (\ref{3.15}), we deduce (\ref{3.12}) and the theorem is proved.
\end{proof}

\section{Further Reverses of the Triangle Inequality}

In this section we point out different additive reverses of the generalised
triangle inequality under simpler conditions for the vectors involved.

The following result holds:

\begin{theorem}
\label{t4.1}Let $\left( H;\left\langle \cdot ,\cdot \right\rangle \right) $
be an inner product space over the real or complex number field $\mathbb{K}$
and $e,x_{i}\in H,$ $i\in \left\{ 1,\dots ,n\right\} $ with $\left\Vert
e\right\Vert =1.$ If $\rho \in \left( 0,1\right) $ and $x_{i},$ $i\in
\left\{ 1,\dots ,n\right\} $ are such that%
\begin{equation}
\left\Vert x_{i}-e\right\Vert \leq \rho \text{ \ for each}\ \ i\in \left\{
1,\dots ,n\right\} ,  \label{4.1}
\end{equation}%
then we have the inequality%
\begin{align}
\left( 0\leq \right) \sum_{i=1}^{n}\left\Vert x_{i}\right\Vert -\left\Vert
\sum_{i=1}^{n}x_{i}\right\Vert & \leq \frac{\rho ^{2}}{\sqrt{1-\rho ^{2}}%
\left( 1+\sqrt{1-\rho ^{2}}\right) }\func{Re}\left\langle
\sum_{i=1}^{n}x_{i},e\right\rangle  \label{4.2} \\
& \left( \leq \frac{\rho ^{2}}{\sqrt{1-\rho ^{2}}\left( 1+\sqrt{1-\rho ^{2}}%
\right) }\left\Vert \sum_{i=1}^{n}x_{i}\right\Vert \right) .  \notag
\end{align}%
The equality holds in (\ref{4.2}) if and only if%
\begin{equation}
\sum_{i=1}^{n}\left\Vert x_{i}\right\Vert \geq \frac{\rho ^{2}}{\sqrt{1-\rho
^{2}}\left( 1+\sqrt{1-\rho ^{2}}\right) }\func{Re}\left\langle
\sum_{i=1}^{n}x_{i},e\right\rangle  \label{4.3}
\end{equation}%
and 
\begin{equation}
\sum_{i=1}^{n}x_{i}=\left( \sum_{i=1}^{n}\left\Vert x_{i}\right\Vert -\frac{%
\rho ^{2}}{\sqrt{1-\rho ^{2}}\left( 1+\sqrt{1-\rho ^{2}}\right) }\func{Re}%
\left\langle \sum_{i=1}^{n}x_{i},e\right\rangle \right) e.  \label{4.4}
\end{equation}
\end{theorem}

\begin{proof}
We know, from the proof of Theorem \ref{t3.1}, that, if (\ref{4.1}) is
fulfilled, then we have the inequality%
\begin{equation*}
\left\Vert x_{i}\right\Vert \leq \frac{1}{\sqrt{1-\rho ^{2}}}\func{Re}%
\left\langle x_{i},e\right\rangle
\end{equation*}%
for each $i\in \left\{ 1,\dots ,n\right\} ,$ implying%
\begin{eqnarray}
\left\Vert x_{i}\right\Vert -\func{Re}\left\langle x_{i},e\right\rangle
&\leq &\left( \frac{1}{\sqrt{1-\rho ^{2}}}-1\right) \func{Re}\left\langle
x_{i},e\right\rangle  \label{4.5} \\
&=&\frac{\rho ^{2}}{\sqrt{1-\rho ^{2}}\left( 1+\sqrt{1-\rho ^{2}}\right) }%
\func{Re}\left\langle x_{i},e\right\rangle  \notag
\end{eqnarray}%
for each $i\in \left\{ 1,\dots ,n\right\} .$

Now, making use of Theorem \ref{t2.1}, for 
\begin{equation*}
k_{i}:=\frac{\rho ^{2}}{\sqrt{1-\rho ^{2}}\left( 1+\sqrt{1-\rho ^{2}}\right) 
}\func{Re}\left\langle x_{i},e\right\rangle ,\ \ \ i\in \left\{ 1,\dots
,n\right\} ,
\end{equation*}%
we easily deduce the conclusion of the theorem.

We omit the details.
\end{proof}

We may state the following result as well:

\begin{theorem}
\label{t4.2}Let $\left( H;\left\langle \cdot ,\cdot \right\rangle \right) $
be an inner product space and $e\in H,$ $M\geq m>0.$ If $x_{i}\in H,$ $i\in
\left\{ 1,\dots ,n\right\} $ are such that either%
\begin{equation}
\func{Re}\left\langle Me-x_{i},x_{i}-me\right\rangle \geq 0,  \label{4.6}
\end{equation}%
or, equivalently,%
\begin{equation}
\left\Vert x_{i}-\frac{M+m}{2}e\right\Vert \leq \frac{1}{2}\left( M-m\right)
\label{4.7}
\end{equation}%
holds for each $i\in \left\{ 1,\dots ,n\right\} ,$ then we have the
inequality%
\begin{align}
\left( 0\leq \right) \sum_{i=1}^{n}\left\Vert x_{i}\right\Vert -\left\Vert
\sum_{i=1}^{n}x_{i}\right\Vert & \leq \frac{\left( \sqrt{M}-\sqrt{m}\right)
^{2}}{2\sqrt{mM}}\func{Re}\left\langle \sum_{i=1}^{n}x_{i},e\right\rangle
\label{4.8} \\
& \left( \leq \frac{\left( \sqrt{M}-\sqrt{m}\right) ^{2}}{2\sqrt{mM}}%
\left\Vert \sum_{i=1}^{n}x_{i}\right\Vert \right) .  \notag
\end{align}%
The equality holds in (\ref{4.8}) if and only if%
\begin{equation}
\sum_{i=1}^{n}\left\Vert x_{i}\right\Vert \geq \frac{\left( \sqrt{M}-\sqrt{m}%
\right) ^{2}}{2\sqrt{mM}}\func{Re}\left\langle
\sum_{i=1}^{n}x_{i},e\right\rangle  \label{4.9}
\end{equation}%
and%
\begin{equation}
\sum_{i=1}^{n}x_{i}=\left( \sum_{i=1}^{n}\left\Vert x_{i}\right\Vert -\frac{%
\left( \sqrt{M}-\sqrt{m}\right) ^{2}}{2\sqrt{mM}}\func{Re}\left\langle
\sum_{i=1}^{n}x_{i},e\right\rangle \right) e.  \label{4.10}
\end{equation}
\end{theorem}

\begin{proof}
We know, from the proof of Theorem \ref{t2.3}, that if (\ref{4.6}) is
fulfilled, then we have the inequality%
\begin{equation*}
\left\Vert x_{i}\right\Vert \leq \frac{M+m}{2\sqrt{mM}}\func{Re}\left\langle
x_{i},e\right\rangle
\end{equation*}%
for each $i\in \left\{ 1,\dots ,n\right\} .$ This is equivalent to%
\begin{equation*}
\left\Vert x_{i}\right\Vert -\func{Re}\left\langle x_{i},e\right\rangle \leq 
\frac{\left( \sqrt{M}-\sqrt{m}\right) ^{2}}{2\sqrt{mM}}\func{Re}\left\langle
x_{i},e\right\rangle
\end{equation*}%
for each $i\in \left\{ 1,\dots ,n\right\} .$

Now, making use of Theorem \ref{t3.1}, we deduce the conclusion of the
theorem. We omit the details.
\end{proof}

\begin{remark}
If one uses Theorem \ref{t3.2} instead of Theorem \ref{t3.1} above, then one
can state the corresponding generalisation for families of orthornormal
vectors of the inequalities (\ref{4.2}) and (\ref{4.8}) respectively. We do
not provide them here.
\end{remark}

Now, on utilising a slightly different approach, we may point out the
following result:

\begin{theorem}
\label{t4.3}Let $\left( H;\left\langle \cdot ,\cdot \right\rangle \right) $
be an inner product space over $\mathbb{K}$ and $e,$ $x_{i}\in H,$ $i\in
\left\{ 1,\dots ,n\right\} $ with $\left\Vert e\right\Vert =1.$ If $r_{i}>0,$
$i\in \left\{ 1,\dots ,n\right\} $ are such that%
\begin{equation}
\left\Vert x_{i}-e\right\Vert \leq r_{i}\text{ \ for each}\ \ i\in \left\{
1,\dots ,n\right\} ,  \label{4.11}
\end{equation}%
then we have the inequality%
\begin{equation}
0\leq \sum_{i=1}^{n}\left\Vert x_{i}\right\Vert -\left\Vert
\sum_{i=1}^{n}x_{i}\right\Vert \leq \frac{1}{2}\sum_{i=1}^{n}r_{i}^{2}.
\label{4.12}
\end{equation}%
The equality holds in (\ref{4.12}) if and only if%
\begin{equation}
\sum_{i=1}^{n}\left\Vert x_{i}\right\Vert \geq \frac{1}{2}%
\sum_{i=1}^{n}r_{i}^{2}  \label{4.13}
\end{equation}%
and%
\begin{equation}
\sum_{i=1}^{n}x_{i}=\left( \sum_{i=1}^{n}\left\Vert x_{i}\right\Vert -\frac{1%
}{2}\sum_{i=1}^{n}r_{i}^{2}\right) e.  \label{4.14}
\end{equation}
\end{theorem}

\begin{proof}
The condition (\ref{4.11}) is clearly equivalent to%
\begin{equation}
\left\Vert x_{i}\right\Vert ^{2}+1\leq \func{Re}\left\langle
x_{i},e\right\rangle +r_{i}^{2}  \label{4.15}
\end{equation}%
for each $i\in \left\{ 1,\dots ,n\right\} .$

Using the elementary inequality%
\begin{equation}
2\left\Vert x_{i}\right\Vert \leq \left\Vert x_{i}\right\Vert ^{2}+1,
\label{4.16}
\end{equation}%
for each $i\in \left\{ 1,\dots ,n\right\} ,$ then, by (\ref{4.15}) and (\ref%
{4.16}), we deduce 
\begin{equation*}
2\left\Vert x_{i}\right\Vert \leq 2\func{Re}\left\langle
x_{i},e\right\rangle +r_{i}^{2},
\end{equation*}%
giving%
\begin{equation}
\left\Vert x_{i}\right\Vert -\func{Re}\left\langle x_{i},e\right\rangle \leq 
\frac{1}{2}r_{i}^{2}  \label{4.17}
\end{equation}%
for each $i\in \left\{ 1,\dots ,n\right\} .$

Now, utilising Theorem \ref{t3.1} for $k_{i}=\frac{1}{2}r_{i}^{2},$ $i\in
\left\{ 1,\dots ,n\right\} ,$ we deduce the desired result. We omit the
details.
\end{proof}

Finally, we may state and prove the following result as well.

\begin{theorem}
\label{t4.4}Let $\left( H;\left\langle \cdot ,\cdot \right\rangle \right) $
be an inner product space over $\mathbb{K}$ and $e,$ $x_{i}\in H,$ $i\in
\left\{ 1,\dots ,n\right\} $ with $\left\Vert e\right\Vert =1.$ If $%
M_{i}\geq m_{i}>0,$ $i\in \left\{ 1,\dots ,n\right\} ,$ are such that%
\begin{equation}
\left\Vert x_{i}-\frac{M_{i}+m_{i}}{2}e\right\Vert \leq \frac{1}{2}\left(
M_{i}-m_{i}\right) ,  \label{4.18}
\end{equation}%
or, equivalently,%
\begin{equation}
\func{Re}\left\langle M_{i}e-x,x-m_{i}e\right\rangle \geq 0  \label{4.19}
\end{equation}%
for each $i\in \left\{ 1,\dots ,n\right\} ,$ then we have the inequality%
\begin{equation}
\left( 0\leq \right) \sum_{i=1}^{n}\left\Vert x_{i}\right\Vert -\left\Vert
\sum_{i=1}^{n}x_{i}\right\Vert \leq \frac{1}{4}\sum_{i=1}^{n}\frac{\left(
M_{i}-m_{i}\right) ^{2}}{M_{i}+m_{i}}.  \label{4.20}
\end{equation}%
The equality holds in (\ref{4.20}) if and only if%
\begin{equation}
\sum_{i=1}^{n}\left\Vert x_{i}\right\Vert \geq \frac{1}{4}\sum_{i=1}^{n}%
\frac{\left( M_{i}-m_{i}\right) ^{2}}{M_{i}+m_{i}}  \label{4.21}
\end{equation}%
and 
\begin{equation}
\sum_{i=1}^{n}x_{i}=\left( \sum_{i=1}^{n}\left\Vert x_{i}\right\Vert -\frac{1%
}{4}\sum_{i=1}^{n}\frac{\left( M_{i}-m_{i}\right) ^{2}}{M_{i}+m_{i}}\right)
e.  \label{4.22}
\end{equation}
\end{theorem}

\begin{proof}
The condition (\ref{4.18}) is equivalent to:%
\begin{equation*}
\left\Vert x_{i}\right\Vert ^{2}+\left( \frac{M_{i}+m_{i}}{2}\right)
^{2}\leq 2\func{Re}\left\langle x_{i},\frac{M_{i}+m_{i}}{2}e\right\rangle +%
\frac{1}{4}\left( M_{i}-m_{i}\right) ^{2}
\end{equation*}%
and since%
\begin{equation*}
2\left( \frac{M_{i}+m_{i}}{2}\right) \left\Vert x_{i}\right\Vert \leq
\left\Vert x_{i}\right\Vert ^{2}+\left( \frac{M_{i}+m_{i}}{2}\right) ^{2},
\end{equation*}%
then we get%
\begin{equation*}
2\left( \frac{M_{i}+m_{i}}{2}\right) \left\Vert x_{i}\right\Vert \leq 2\cdot 
\frac{M_{i}+m_{i}}{2}\func{Re}\left\langle x_{i},e\right\rangle +\frac{1}{4}%
\left( M_{i}-m_{i}\right) ^{2},
\end{equation*}%
or, equivalently,%
\begin{equation*}
\left\Vert x_{i}\right\Vert -\func{Re}\left\langle x_{i},e\right\rangle \leq 
\frac{1}{4}\cdot \frac{\left( M_{i}-m_{i}\right) ^{2}}{M_{i}+m_{i}}
\end{equation*}%
for each $i\in \left\{ 1,\dots ,n\right\} .$

Now, making use of Theorem \ref{t3.1} for $k_{i}:=\frac{1}{4}\cdot \frac{%
\left( M_{i}-m_{i}\right) ^{2}}{M_{i}+m_{i}},$ $i\in \left\{ 1,\dots
,n\right\} ,$ we deduce the desired result.
\end{proof}

\begin{remark}
If one uses Theorem \ref{t3.2} instead of Theorem \ref{t3.1} above, then one
can state the corresponding generalisation for families of orthornormal
vectors of the inequalities in (\ref{4.12}) and (\ref{4.20}) respectively.
We omit the details.
\end{remark}

\section{Reverses of Schwarz Inequality}

In this section we outline a procedure showing how some of the above results
for triangle inequality may be employed to obtain reverses for the
celebrated Schwarz inequality.

For $a\in H,$ $\left\Vert a\right\Vert =1$ and $r\in \left( 0,1\right) $
define the closed ball 
\begin{equation*}
\overline{D}\left( a,r\right) :=\left\{ x\in H,\left\Vert x-a\right\Vert
\leq r\right\} .
\end{equation*}%
The following reverse of the Schwarz inequality holds:

\begin{proposition}
\label{p5.1} If $x,y\in \overline{D}\left( a,r\right) $ with $a\in H,$ $%
\left\Vert a\right\Vert =1$ and $r\in \left( 0,1\right) ,$ then we have the
inequality%
\begin{equation}
\left( 0\leq \right) \frac{\left\Vert x\right\Vert \left\Vert y\right\Vert -%
\func{Re}\left\langle x,y\right\rangle }{\left( \left\Vert x\right\Vert
+\left\Vert y\right\Vert \right) ^{2}}\leq \frac{1}{2}r^{2}.  \label{5.1}
\end{equation}%
The constant $\frac{1}{2}$ in $\left( \ref{5.1}\right) $ is best possible in
the sense that it cannot be replaced by a smaller quantity.
\end{proposition}

\begin{proof}
Using Theorem \ref{t2.1} for $x_{1}=x,x_{2}=y,\rho =r,$ we have 
\begin{equation}
\sqrt{1-r^{2}}\left( \left\Vert x\right\Vert +\left\Vert y\right\Vert
\right) \leq \left\Vert x+y\right\Vert .  \label{5.2.1}
\end{equation}%
Taking the square in (\ref{5.2.1}) we deduce 
\begin{equation*}
\left( 1-r^{2}\right) \left( \left\Vert x\right\Vert ^{2}+2\left\Vert
x\right\Vert \left\Vert y\right\Vert +\left\Vert y\right\Vert ^{2}\right)
\leq \left\Vert x\right\Vert ^{2}+2\func{Re}\left\langle x,y\right\rangle
+\left\Vert y\right\Vert ^{2}
\end{equation*}%
which is clearly equivalent to (\ref{5.1}).

Now, assume that (\ref{5.1}) holds with a constant $C>0$ instead of $\frac{1%
}{2},i.e.,$%
\begin{equation}
\frac{\left\Vert x\right\Vert \left\Vert y\right\Vert -\func{Re}\left\langle
x,y\right\rangle }{\left( \left\Vert x\right\Vert +\left\Vert y\right\Vert
\right) ^{2}}\leq Cr^{2}  \label{5.3}
\end{equation}%
provided $x,y\in \overline{D}\left( a,r\right) $ with $a\in H,$ $\left\Vert
a\right\Vert =1$ and $r\in \left( 0,1\right) .$

Let $e\in H$ with $\left\Vert e\right\Vert =1$ and $e\perp a.$ Define $%
x=a+re,y=a-re.$ Then 
\begin{equation*}
\left\Vert x\right\Vert =\sqrt{1+r^{2}}=\left\Vert y\right\Vert ,\text{ }%
\func{Re}\left\langle x,y\right\rangle =1-r^{2}
\end{equation*}%
and thus, from (\ref{5.3}), we have%
\begin{equation*}
\frac{1+r^{2}-\left( 1-r^{2}\right) }{\left( 2\sqrt{1+r^{2}}\right) ^{2}}%
\leq Cr^{2}
\end{equation*}%
giving 
\begin{equation*}
\frac{1}{2}\leq \left( 1+r^{2}\right) C
\end{equation*}%
for any $r\in \left( 0,1\right) .$ If in this inequality we let $%
r\rightarrow 0+,$ then we get $C\geq \frac{1}{2}$ and the proposition is
proved.
\end{proof}

In a similar way, by the use of Theorem \ref{t2.3}, we may prove the
following reverse of the Schwarz inequality as well:

\begin{proposition}
\label{p5.2} If $a\in H,$ $\left\Vert a\right\Vert =1,$ $M\geq m>0$ and $%
x,y\in H$ are so that either%
\begin{equation*}
\func{Re}\left\langle Ma-x,x-ma\right\rangle ,\func{Re}\left\langle
Ma-y,y-ma\right\rangle \geq 0
\end{equation*}%
or, equivalently,%
\begin{equation*}
\left\Vert x-\frac{m+M}{2}a\right\Vert ,\left\Vert y-\frac{m+M}{2}%
a\right\Vert \leq \frac{1}{2}\left( M-m\right)
\end{equation*}%
hold, then 
\begin{equation*}
\left( 0\leq \right) \frac{\left\Vert x\right\Vert \left\Vert y\right\Vert -%
\func{Re}\left\langle x,y\right\rangle }{\left( \left\Vert x\right\Vert
+\left\Vert y\right\Vert \right) ^{2}}\leq \frac{1}{2}\left( \frac{M-m}{M+m}%
\right) ^{2}.
\end{equation*}%
The constant $\frac{1}{2}$ cannot be replaced by a smaller quantity.
\end{proposition}

\begin{remark}
On utilising Theorem \ref{t2.2} and Theorem \ref{t2.4}, we may deduce some
similar reverses of Schwarz inequality provided $x,y\in \cap _{k=1}^{m}%
\overline{D}\left( a_{k},\rho _{k}\right) ,$ assumed not to be empty, where $%
a_{1},...,a_{n}$ are orthonormal vectors in $H$ and $\rho _{k}\in \left(
0,1\right) $ for $k\in \left\{ 1,...,m\right\} .$ We omit the details.
\end{remark}

\begin{remark}
For various different reverses of Schwarz inequality in inner product
spaces, see the recent survey \cite{SSD}, that is available as a preprint in
Mathematical Ar$\chi $iv, where further references are given.
\end{remark}

\section{Applications for Vector-Valued Integral Inequalities}

Let $\left( H;\left\langle \cdot ,\cdot \right\rangle \right) $ be a Hilbert
space over the real or complex number field, $\left[ a,b\right] $ a compact
interval in $\mathbb{R}$ and $\eta :\left[ a,b\right] \rightarrow \lbrack
0,\infty )$ a Lebesgue integrable function on $\left[ a,b\right] $ with the
property that $\int_{a}^{b}\eta \left( t\right) dt=1.$ If, by $L_{\eta
}\left( \left[ a,b\right] ;H\right) $ we denote the Hilbert space of all
Bochner measurable functions $f:\left[ a,b\right] \rightarrow H$ with the
property that $\int_{a}^{b}\eta \left( t\right) \left\Vert f\left( t\right)
\right\Vert ^{2}dt<\infty ,$ then the norm $\left\Vert \cdot \right\Vert
_{\eta }$ of this space is generated by the inner product $\left\langle
\cdot ,\cdot \right\rangle _{\eta }:H\times H\rightarrow \mathbb{K}$ defined
by 
\begin{equation*}
\left\langle f,g\right\rangle _{\eta }:=\int_{a}^{b}\eta \left( t\right)
\left\langle f\left( t\right) ,g\left( t\right) \right\rangle dt.
\end{equation*}%
The following proposition providing a reverse of the integral generalised
triangle inequality may be stated.

\begin{proposition}
\label{p6.1} Let $\left( H;\left\langle \cdot ,\cdot \right\rangle \right) $
be a Hilbert space and $\eta :\left[ a,b\right] \rightarrow \lbrack 0,\infty
)$ as above. If $g\in L_{\eta }\left( \left[ a,b\right] ;H\right) $ is so
that $\int_{a}^{b}\eta \left( t\right) \left\Vert g\left( t\right)
\right\Vert ^{2}dt=1$ and $f_{i}\in L_{\eta }\left( \left[ a,b\right]
;H\right) ,i\in \left\{ 1,\dots ,n\right\} ,$ $\rho \in \left( 0,1\right) $
are so that%
\begin{equation}
\left\Vert f_{i}\left( t\right) -g\left( t\right) \right\Vert \leq \rho
\label{6.1}
\end{equation}%
for a.e. $t\in \left[ a,b\right] $ and each $i\in \left\{ 1,\dots ,n\right\}
,$ then we have the inequality%
\begin{equation}
\sqrt{1-\rho ^{2}}\sum_{i=1}^{n}\left( \int_{a}^{b}\eta \left( t\right)
\left\Vert f_{i}\left( t\right) \right\Vert ^{2}dt\right) ^{1/2}\leq \left(
\int_{a}^{b}\eta \left( t\right) \left\Vert \sum_{i=1}^{n}f_{i}\left(
t\right) \right\Vert ^{2}dt\right) ^{1/2}.  \label{6.2}
\end{equation}%
The case of equality holds in (\ref{6.2}) if and only if%
\begin{equation*}
\sum_{i=1}^{n}f_{i}\left( t\right) =\sqrt{1-\rho ^{2}}\sum_{i=1}^{n}\left(
\int_{a}^{b}\eta \left( t\right) \left\Vert f_{i}\left( t\right) \right\Vert
^{2}dt\right) ^{1/2}\cdot g\left( t\right)
\end{equation*}%
for a.e. $t\in \left[ a,b\right] .$
\end{proposition}

\begin{proof}
Observe, by (\ref{6.2}), that%
\begin{eqnarray*}
\left\Vert f_{i}-g\right\Vert _{\eta } &=&\left( \int_{a}^{b}\eta \left(
t\right) \left\Vert f_{i}\left( t\right) -g\left( t\right) \right\Vert
^{2}dt\right) ^{1/2} \\
&\leq &\left( \int_{a}^{b}\eta \left( t\right) \rho ^{2}dt\right) ^{1/2}=\rho
\end{eqnarray*}%
for each $i\in \left\{ 1,\dots ,n\right\} .$ Applying Theorem \ref{t2.1} for
the Hilbert space $L_{\eta }\left( \left[ a,b\right] ;H\right) ,$ we deduce
the desired result.
\end{proof}

The following result may be stated as well.

\begin{proposition}
\label{p6.2} Let $H,\eta ,g$ be as in Proposition \ref{p6.1}. If $f_{i}\in
L_{\eta }\left( \left[ a,b\right] ;H\right) ,i\in \left\{ 1,\dots ,n\right\} 
$ and $M\geq m>0$ are so that either%
\begin{equation*}
\func{Re}\left\langle Mg\left( t\right) -f_{i}\left( t\right) ,f_{i}\left(
t\right) -mg\left( t\right) \right\rangle \geq 0
\end{equation*}%
or, equivalently, 
\begin{equation*}
\left\Vert f_{i}\left( t\right) -\frac{m+M}{2}g\left( t\right) \right\Vert
\leq \frac{1}{2}\left( M-m\right)
\end{equation*}%
for a.e. $t\in \left[ a,b\right] $ and each $i\in \left\{ 1,\dots ,n\right\}
,$ then we have the inequality%
\begin{equation}
\frac{2\sqrt{mM}}{m+M}\sum_{i=1}^{n}\left( \int_{a}^{b}\eta \left( t\right)
\left\Vert f_{i}\left( t\right) \right\Vert ^{2}dt\right) ^{1/2}\leq \left(
\int_{a}^{b}\eta \left( t\right) \left\Vert \sum_{i=1}^{n}f_{i}\left(
t\right) \right\Vert ^{2}dt\right) ^{1/2}.  \label{6.3}
\end{equation}%
The equality holds in (\ref{6.3}) if and only if%
\begin{equation*}
\sum_{i=1}^{n}f_{i}\left( t\right) =\frac{2\sqrt{mM}}{m+M}%
\sum_{i=1}^{n}\left( \int_{a}^{b}\eta \left( t\right) \left\Vert f_{i}\left(
t\right) \right\Vert ^{2}dt\right) ^{1/2}\cdot g\left( t\right) ,
\end{equation*}%
for a.e. $t\in \left[ a,b\right] .$
\end{proposition}

\begin{remark}
Similar integral inequalities may be stated on utilising the inequalities
for inner products and norms obtained above, but we do not mention them here.
\end{remark}

\end{document}